\magnification=\magstep1

\font\TITLE=cmbx10 scaled \magstep1 \font\AUT=cmcsc10
\font\SUBT=cmssdc10 
\font\aop=cmr9 \font\eop=cmr8 \overfullrule=0pt

\def\n{\noindent} \def\br{\buildrel}
\def\bs{\bigskip} \def\ms{\medskip} \def\sm{\smallskip} \def\ol{\overline}
\def\cA{{\cal A}}  \def\cF{{\cal F}} \def\cK{{\cal K}}
\def\cL{{\cal L}}
\def\cN{{\cal N}}
\def\la{\lambda} \def\de{\delta} \def\ve{\varepsilon} \def\QED{\hfill{\eop
QED}} 
\def\RA{\,\Rightarrow\,} 

\vglue 15mm
\centerline {\TITLE Self-Induced Compactness in Banach Spaces} \bs\ms
\centerline {\AUT P.G.Casazza and H.Jarchow} 

\vskip 10mm
\centerline {\SUBT Introduction}
\bs\sm
The question which led to the title of this note is the following: 
\ms
{\it If $X$ is a Banach space and $K$ is a compact subset of $X$, is it
possible to find a 
compact, or even approximable, operator $v:X\to X$ such that
$K\subset\ol{v(B_X)}$?}
\ms
This question was first posed by P.G.Dixon [6] in connection with
investigating the problem of
the existence of approximate identities in certain operator algebras. We
shall provide a couple 
of observations related to the above question and give in particular a
negative answer in case of
approximable operators.
\sm
We shall also provide the first examples of Banach spaces having the
approximation property but 
failing the bounded compact approximation property though all of their
duals do even have the 
metric compact approximation property. 

\vglue 10mm\n
\centerline {\SUBT Approximate Identities} \bs\sm

A {\it left approximate identity} (LAI) in a Banach algebra $A$ is a net
$(e_i)_{i\in I}$ in $A$
such that $lim_{i\in I}\|e_ix-x\|=0$ for each $x\in A$. If there is a
$\la>0$ such that in
addition $\|e_i\|\le\la$, then $(e_i)_{i\in I}$ is called a $\la$-{\it
bounded LAI\/} 
($\la$-BLAI); note that necessarily $\la\ge 1$. We say that $(e_i)_{i\in
I}$ is a BLAI
if it is a $\la$-BLAI for some $\la$.
\sm
{\it Right approximate identities\/} (RAI's) and ($\la$-){\it bounded\/}
RAI's (($\la$-)BRAI's) 
are defined analogously.
\sm
A result which, when seen in the context of operator algebras, is of
particular interest for 
Banach space theory, is P.J.Cohen's [4] factorization theorem: 
\sm\n
{\it If the Banach algebra $A$ has a BLAI, then every $z\in A$ can be
written as a product $z=xy$
with $x,y\in A$; in addition, $y$ exists in the closed left ideal generated
by $z$ and can be 
chosen so that $\|z-y\|\le\de$, where $\de$ was given previously.} 
\sm
See the book [2] of F.F.Bonsall and J.Duncan for details. 
\ms
An interesting case occurs when $A$ is of the form $\cA(X)$ where $\cA$ is
any Banach operator 
ideal in the sense of A.Pietsch [13] and $X$ is a Banach space. However, in
such generality 
approximate identities have hardly been discussed; attention has been
focused on the cases where 
$\cA$ is either
$$\cK~,$$
the ideal of all compact operators, or
$$\ol\cF~,$$

\goodbreak\n
the ideal of all {\it approximable operators,} i.e., uniform limits of
operators of finite rank
operators. Of course,
$$\cF$$
will be used to denote the ideal of all finite rank operators between
Banach spaces. Dixon's 
paper [6] is dealing with precisely these ideals, and one of his main
results is the following:
\ms\n
\proclaim Theorem A. A Banach space $X$ has the $\la$-BAP [resp. the
$\la$-BCAP] if and only if
$\ol\cF(X)$ [resp. $\cK(X)$] has a $\la$-BLAI. 

\ms
Recall that a Banach space $X$ has {\it AP (``approximation property")\/}
[resp.{\it CAP 
(``compact approximation property")\/}] if, for any compact subset $K$ of
$X$ and any $\ve>0$, 
there is an operator $u\in\cF(X)$ [resp. $u\in\cK(X)$] such that
$\|ux-x\|<\ve$ for all $x\in K$. 
If there is a $\la(\ge 1)$ such that we always can arrange for
$\|u\|\le\la$, then we say that 
$X$ has $\la$-{\it BAP\/} [resp. $\la$-{\it BCAP\/}], with ``B" being
shorthand for ``bounded", 
of course. Usually, {\it BAP\/} [resp. {\it BCAP\/}] is used when we only
know that we have 
$\la$-BAP [resp. $\la$-BCAP] for some $\la$. The case $\la=1$ corresponds to
what is usually
called {\it MAP (``metric approximation property")\/} and {\it MCAP
(``metric compact 
approximation property"),} respectively.
\ms
Only recently C.Samuel [16] and N.Gr\o nb\ae k and G.A.Willis [7] addressed
the corresponding
problem of existence of RAI's. One of the main results in [7] is the
following companion result 
of the preceding theorem:
\ms
\proclaim Theorem B. Let $X$ be a Banach space. Then $X^*$ has the $\la$-BAP if
and only if ${\ol\cF}(X)$ has a $\la$-BRAI.

\ms
The situation is less pleasant when $\ol\cF$ is replaced by $\cK.$ It is
still true that if 
$\cK(X)$ has a $\la$-BRAI, then $X$ has $\la$-BCAP. However, the converse
fails. As we shall see, 
there is a Banach space $X$ with AP, but which fails BCAP, while $X^*,
X^{**},\dots$ are all 
separable and have MCAP. If $\ol\cF(X)=\cK(X)$ has a BRAI, then $X^*$ would
have BAP, and so $X$ 
should have BAP as well -- but it doesn't.
\ms
To construct our example we need some results from the literature. Recall
that a Banach space $X$ 
has {\it shrinking\/} $\la$-CAP if there is a net $u_\alpha\in\cK(X)$,
$\|u_\alpha\|\le\la$,
$(u_{\alpha})$ converges strongly to the identity on $X\; {\underline
{and}}\; (u^*_{\alpha})$
converges strongly to the identity on $X^*$. In his memoir [8], A.
Grothendieck showed that a
separable dual space with AP has MAP. The corresponding question for CAP is
still open. However,
an alternative proof of Grothendieck's result due to Lindenstrauss and
Tzafriri [11] will work 
for shrinking CAP (see Cho and Johnson [3]). 
\ms\n
\proclaim Theorem C. If $X^*$ is separable and has CAP given by
$w^*$-continuous operators, then 
$X$ and $X^*$ have MCAP.

\ms
It follows that,
\ms
\proclaim Corollary D. If $X^*$ is separable and $X$ has shrinking CAP,
then for every equivalent
norm $|\cdot|$ on $X$, both $X$ and $(X,|\cdot|)^*$ have MCAP. 

\ms
We are now ready for:

\proclaim Theorem 1. There is a Banach space $X$ which has AP, but fails
BCAP, while $X^*, X^{**},
\dots$ are all separable and have MCAP.

\goodbreak\n
{\sl Proof:}~~~Let $Y$ be a separable reflexive Banach space with CAP which
fails to have AP, see
G.Willis [19]. Choose a Banach space $Z$ so that $Z^{**}$ has a basis and
$Z^{**}/Z\approx Y$. 
(This result and the others used here can be found in [11]). Now,
$Z^{***}\approx Z^*\oplus Y^*$ 
fails AP but has shrinking CAP since (1)~~$Z^{*}$ has a shrinking basis and
(2)~~$Y$ is reflexive 
and has CAP. It follows that $Y$ has shrinking CAP. By a construction of
Figiel and Johnson
([11], p.42), for each $n$, there is an equivalent norm $|\cdot|_n$ on
$Z^{**}$ so that 
$(Z^{**},|\cdot|_n)$ fails $n-$BAP. But $(Z^{**},|\cdot|_n)^*$ still has
shrinking CAP (being 
isomorphic to $Z^{***}$). So by Corollary D, this space and its dual have
MCAP. (A similar 
argument shows, in fact, that all of its duals have MCAP.)
\sm
Let
$$X\,=\,\left(\sum^{\infty}_{n=1}\oplus(Z^{**},|\cdot|_n)\right)_{\ell_{2}}~.$$
Then $X$ has AP and fails BAP. So $X$ fails BCAP. (It is easily seen that
$X$ has $\la$-BAP if 
and only if $X$ has AP and $\la$-BCAP). The spaces $X^*, X^{**},\dots$ are
all $\ell_2$-sums of 
Banach spaces having MCAP, and hence have MCAP. This completes the
construction.\QED

\vglue 10mm\n
\centerline {\SUBT The Properties (${\ol\cF}$) and ($\cK$)} \bs\sm

We shall now concentrate on left approximate units. There is a natural
question related to 
Theorem A: if the same $\la$ works for the BLAI and for the B(C)AP,
couldn't it be simultaneously 
`eliminated' on both sides ? More precisely, is it true that ${\cal \ol
F}(X)$ resp. $\cK(X)$ has 
a LAI if and only if $X$ has AP resp. CAP ?
\sm
Let $\cA$ be ${\cal \ol F}$ or $\cK$. We say that a Banach space $X$ has
the {\it property 
$(\cA)$} if for each compact subset $K$ of $X$ there is an operator
$u\in\cA(X)$ such that 
$K\subset\ol{u(B_X)}$. This is what ``self-induced compactness" in the
title is  referring to. Of 
course, the concept can be generalized in many directions, but we prefer to
stay with the present 
setup.
\sm
The following results are again due to Dixon [6]: 
\ms
\proclaim Theorem E. Let $X$ be a Banach space. 
\sm\n
\item{\rm(a)~}If $X$ has AP [CAP], then $\ol\cF(X)$ [$\cK(X)$] has a LAI.
\sm\n
\item{\rm(b)~}If $X$ has property ($\ol\cF$) [($\cK$)] and if $\ol\cF(X)$
[$\cK(X)$] has a LAI, then $X$ has AP [CAP].
\sm\n
\item{\rm(c)~}If $X$ has BAP [BCAP], then it has property ($\ol\cF$) [$(\cK)$].

\ms
All we need to get started is a workable condition which is equivalent to
property ($\cA$). This 
is elementary:
\ms\n
\proclaim Proposition 1. The following statements about the Banach space
$X$ are equivalent:
\sm
\item{\rm(i)~} $X$ has property ($\cal A$). 
\sm
\item{\rm(ii)~} Given $u\in\cK(\ell_1,X)$, there are operators $v\in\cA(X)$
and $w_n\in\cL(\ell_1,X)$
$(n\in{\bf N})$ such that $\lim_{n\to\infty}\|u-vw_n\|=0$. 
\sm
\item{\rm(iii)~} same as (ii), but $w_n\in\cK(\ell_1,X)$. 
\sm
\item{\rm(iv)~}For each $u\in\cK(\ell_1,X)$, there are operators
$v\in\cA(X)$ and $w\in
\cL(\ell_1,X^{**})$ such that $u=v^{**}w.$
\sm
\item{\rm(v)~} same as (iv), but $w\in\cK(\ell_1,X^{**})$. 

\goodbreak
Since we are working with (weakly) compact operators $v:X\to X$, we may and
shall consider
$v^{**}$ as an operator $X^{**}\to X$; accordingly we have 
$$\ol{v(B_X)}=v^{**}(B_{X^{**}}).$$
\ms\n
{\sl Proof.}~~~(i)$\RA$(ii):~~~Given $u\in\cK(\ell_1,X)$, we can find for
each $n\in {\bf N}$ 
vectors $x_k\in B_X$ such that $\|vx_k-ue_k\|\le n^{-1}$ for all $k\in{\bf
N}$, and then define 
$w_n:\ell_1\to X$ via $w_ne_k=x_k$ for each $k$.
\sm\n
(ii)$\RA$(iii):~~~Any $u\in\cK(\ell_1,X)$ can be written as $u=u_1u_2$
where $u_1:\ell_1\to X$ 
and $u_2:\ell_1\to \ell_1$ are compact operators. 
\sm\n
(iii)$\RA$(iv):~~~We may assume that $\|w_n\|\le 1$ for each $n$. Let 
$\cal U$ be a free 
ultrafilter on $\bf N$ and define $w:\ell_1\to X^{**}$ by
$w(\xi)=\lim_{\cal U}w_n\xi$, the 
limit being taken in the weak*topology of $X^{**}$. This is a bounded
linear operator, and 
$u=v^{**}w.$ 
\sm\n
(iv)$\RA$(v) is obtained as was (ii)$\RA$(iii), so we are left  with
(v)$\RA$(i): Let 
$K\subset X$ be compact. Then $K\subset\ol{\rm conv}\,\{x_n:n\in\bf N\}$
for some null sequence 
$(x_n)$ in $K$. The operator $u:\ell_1\to X$ defined by
$ue_n:=x_n$ for each $n$ is compact, so there are $v\in\cA(X)$ and
$w\in\cK(\ell_1,X^{**})$ such
that $u=v^{**}w$. We may assume that $\|w\|\le 1$ so that
$K\subset\ol{u(B_{\ell_1})}\subset 
v^{**}(B_{X^{**}})$.\QED
\ms
Since $\ell_1$ enjoys the lifting property, a Banach space $X$ has property
$(\cA)$ whenever the
following applies: no matter how we choose the Banach space $Y$ and the
operator $u\in\cA(Y,X)$, 
we can find a quotient $Q$ of $X$ along with operators $v\in\cA(Q,X)$ and
$w\in\cL(Y,Q^{**})$ 
such that $u=v^{**}w$. It is interesting to note that in case $\cA=\cK$,
the preceding 
proposition allows to extend this almost to a characterization. We have the
following weak 
version of the Cohen factorization theorem.
\ms
\proclaim Proposition 2. Let $X$ be a Banach space. 
\sm\n
\item{\rm(a)~}If $X$ has property $(\cK)$ then, given any Banach space $Y$,
every operator $u\in\cK(Y,X)$ admits a compact factorization through some
quotient $Q$ of
$X^{**}$: there are operators $v\in\cK(Q,X)$ and $w\in\cK(Y,Q)$ such that
$u=vw$. 
\sm\n
\item{\rm(b)~}Suppose there is, for every Banach space $Y$ and every
$u\in\cK(Y,X)$, a quotient
space $Q$ of $X$ together with operators $v\in\cK(Q,X)$ and
$w\in\cK(Y,Q^{**})$ such that
$u=v^{**}w$. Then $X$ has property $(\cK)$. 

\ms\n
{\sl Proof.}~~~(a)~~~Suppose that $X$ has $(\cK)$ and let $Y$ and
$u\in\cK(Y,X)$ be given. Of 
course we may assume that $Y$ is sparable; so we can work with a quotient
map $q:\ell_1\to Y.$ 
Thanks to Proposition 1 there are operators $v_0\in\cK(X)$ and
$w_0\in\cK(\ell_1,X^{**})$ such 
that $uq=v_0^{**}w_0.$ Set $Q=X^{**}/ker(v_0^{**})$, let $p:X^{**}\to Q$ be
the quotient map, and 
let
$v\in\cK(Q,X)$ be such that $vp=v_0^{**}.$ As $v$ is injective, there is a
$w\in\cK(Y,Q)$ such 
that $wq=pw_0$. Note that $u=vw.$
\sm\n
(b)~~~Apply the hypothesis to any $\cK(\ell_1,X)$ and use the lifting
property of $\ell_1$.\QED
\ms
We do not have a corresponding result for the property $(\ol\cF)$. \ms
We continue by giving a number of immediate consequences of Proposition 1.
\sm
If the Banach space $X$ admits a quotient which is isomorphic to $\ell_1$
then, by the lifting 
property of $\ell_1$, this quotient is isomorphic to a complemented
subspace of $X$. So we may
state:
\bs
\proclaim Corollary 1. Any Banach space which admits a quotient isomorphic
to $\ell_1$ has the property ($\ol\cF).$

\ms
If $X$ fails CAP then $X\oplus\ell_1$ fails CAP and has $\ell_1$ as a
quotient space. It follows
from Theorem E that for such a space $X$ neither $\cK(X)$ nor $\ol\cF(X)$
can have a LAI. We can easily extend the list of such examples. \ms
\proclaim Lemma. Let $X,\,Y$ and $Z$ be Banach spaces. Suppose that $X$ is
a quotient of $Y$, that
$Z$ is a quotient of $X$, and that every $u\in\cK(\ell_1,Y)$ can be written
$u=vw$ where
$w\in\cK(\ell_1,Z)$ and $v\in\cA(Z,X)$. Then X has the property $(\cA)$.

\ms
This is an immediately consequence of the fact that $\ell_1$ has the
compact lifting property.
\ms
\proclaim Corollary 2. Let $X$ be a subspace of $c_0$ or $\ell_p$ $(1\le
p<\infty)$. Then $X^*$ has property $(\ol\cF).$

\ms
But $X^*$ may well fail to have CAP in which case again neither $\ol\cF(X)$
nor $\cK(X)$ can have 
a LAI.
\sm
The proof is immediate from the lemma and the fact that every infinite
dimensional subspace of 
$c_0$ or $\ell_p$ contains a subspace which is isomorphic to $c_0$ resp.
$\ell_p$ and 
complemented in the whole space.
\ms
We do not know if every subspace of $c_0$ or $\ell_p$ ($1\le p<\infty$) has
$(\cK)$ or 
$(\ol\cF)$. Also, we do not know how $(\cA)$ behaves with respect to duality. 
\sm
W.B.Johnson has proved in [10] that there is a separable Banach space,
$C_1$, such that if $Z$ is
any separable Banach space, then $Z^*$ is isometrically isomorphic to a
norm-one complemented 
subspace of ${C_1}^*$. It follows that ${C_1}^*$ fails CAP but it has the
property $(\ol\cF)$ 
since $\ell_1$ is complemented in ${C_1}^*$. Consequently, neither
$\cK({C_1}^*)$ nor 
$\ol\cF({C_1}^*)$ can have a LAI.
\sm
On the basis of this it is tempting to conjecture that the property $(\cal
A)$ is preserved under
the formation of complemented subspaces. The above observation about duals
of subspaces of $c_0$ 
and $\ell_p$ would then appear as a consequence of ${C_1}^*$'s property of
having $(\ol\cF)$; in 
fact the dual of any separable Banach space would have $(\ol\cF)$. We shall
now see that this is 
not the case.
\sm
Let $\Gamma_2$ denote the ideal of all Banach space operators which factor
through some Hilbert
space.
\ms
\proclaim Proposition 3. Suppose that the Banach space $X$ has the property
($\cal A$). If
$\cA(X)\subset {\Gamma_2}(X)$, then $X$ is isomorphic to a Hilbert space.

\ms\n
{\sl Proof.}~~~In fact, Grothendieck's Inequality informs us that
$\cK(\ell_1,X)$ consists of
absolutely summing operators only. By trace duality, however, this can only
happen when $X$ is
isomorphic to a Hilbert space.\QED
\ms
It was shown by G.Pisier [14] that ${\ol\cF}(X)$ is contained in
$\Gamma_2(X)$ whenever $X$ and 
$X^*$ both have cotype 2, and a few years later, he proved [15] that every
cotype $2$ space $Z$ 
embeds into a non-hilbertian space $P_Z$ which, together with its dual, has
cotype 2. Such a 
space necessarily fails AP. Thus:
\ms
\proclaim Corollary 3. The Pisier spaces $P_Z$ fail to have property $(\ol\cF)$.

\ms
The same is true for all the duals of the spaces $P_Z$. \ms
But thanks to Corollary 1, the spaces $P_Z\oplus\ell_1$ enjoy property
$(\ol\cF)$. Conclusion:
\ms
\proclaim Corollary 4. The property $(\ol\cF)$ is not preserved when
passing to complemented
subspaces.

\ms
Another way to obtain this is by using Johnson's universal space ${C_1}^*$:
if $Z$ is separable,
then Pisier's construction leads to a separable space $P_Z$, so that
${P_Z}^*$ is complemented in
${C_1}^*$.

\ms
Actually, the spaces $P_Z$ enjoy even more exotic properties. For example,
${\ol\cF}(P_Z)$ 
coincides with ${\cal N}(P_Z)$, the algebra of
all {\it nuclear} operators $u:P_Z\to P_Z$. We do 
not know if ${\ol\cF}(P_Z)=\cN(P_Z)$ can have a LAI; actually, the question
of what the
meaning of the existence of LAI's in $\cN(X)$ is in terms of $X$ doesn't
seem to have been 
investigated. However, it was recently shown by Selivanov [17] that
$\cN(X)$ has a BLAI if and 
only if $X$ is finite dimensional; see also Dales and Jarchow in [5].
\sm
The arguments employed before also show that $\Gamma_2(P_Z)$ doesn't have a LAI.

\vglue 10mm\n
\centerline {\SUBT Odds and Ends}
\bs\sm

We start by listing a few open problems: 
\ms\n
(a) Are there Banach spaces failing $(\cK)$? 
\sm\n
(b) Are there Banach spaces failing $(\cal K)$ such that $\cK(X)$ has
(doesn't have) a LAI?
\sm\n
(c) Can any of the algebras ${\ol\cF}(P_Z)$ have a LAI? 
\sm\n
(d) When do the algebras ${\cN}(X)$, $\Gamma_2(X)$,..... have a LAI? 
\sm\n
(e) How do the properties $(\cK)$ and $(\ol\cF)$ behave with respect to duality?
\sm\n
(f) Regarding $(\cK)$ and $(\ol\cF)$, what can be said when the underlying
Banach space X is a Banach lattice, $H^{\infty}$, a $C^*$-algebra, .... ? 
\ms
Here are some further ideas which lead to many more problems. 
\sm
The property $(\cA)$ can be generalized as follows. Given a Banach space
$X$, let
$$\cA_X$$
be the collection of all Banach spaces $Z$ such that for each compact
subset $K$ of $X$ there is
an operator $u\in\cA(Z,X)$ such that $K\subset {\ol{u(B_Z)}}$; again $\cA$
is $\cK$ or $\ol\cF$.
\ms
It seems plausible that investigation of such a concept could be helpful in
understanding 
compactness in general Banach spaces through known characterizations of
compactness in
e.g. classical spaces. Not much, however, is known, and what is known
indicates that the picture 
will by no means be easy to understand.
\sm
There are Banach spaces $Z$ which belong to $\bigcap_X{\ol\cF_X}$: $\ell_1$
and the duals of Johnson's spaces $C_p$ provide examples. 
\ms
If $X$ has BAP [resp. BCAP], then $X$ belongs to $\ol\cF_X$ [resp.
$\cK_X$], whereas the Pisier
spaces $P_Z$ satisfy $P_Z\notin\ol\cF_{P_Z}$. 
\ms
These spaces also satisfy $P_Z\notin\cK_{{P_Z}^*}$. In fact, K.John [9] has
shown that every 
compact operator $P_Z\to{P_Z}^*$ is nuclear, so that $P_Z\in\cK_{{P_Z}^*}$
would entail 
$\cK(\ell_1,{P_Z}^*)={\cN}(\ell_1,{P_Z}^*)$ which cannot be reconciled with
$P_Z$ being 
infinite-dimensional.
\ms
On the other hand, there are Banach spaces $X$ such that
$X\in\ol\cF_{X^*}$: think of $X=\ell_1$, 
$X=C_p$, $X=\ell_2$,.... Can one characterize such spaces? 
\ms
By the same type of argument we get that if $X$ is a non-hilbertian cotype
2 space, then 
$\ol\cF_X$ cannot contain any $Z$ such that $Z^*$ has cotype 2. In fact,
otherwise we
would get $\cK(\ell_1,X)=\Gamma_2(\ell_1,X)$ (cf. [14]), and this is only
possible if $X$ is 
isomorphic to a Hilbert space. Similar, if $X$ has cotype 2 and if $\cK_X$
contains a Banach 
space of type 2, then $X$ must be hilbertian.
\ms
All these example revolve around Hilbert space and amount to ending up with
the conclusion that
$\cK_X$ contains a Hilbert space iff $X$ is isomorphic to Hilbert space.
\ms
We claim that, if $X$ is hilbertian, then $\cK_X$ contains even all dual
Banach spaces. (All 
spaces are supposed to be infinite dimensional). This can be seen as follows:
\ms
Let $u:\ell_1\to X$ be a compact operator. Since $X$ is hilbertian,
$u:\ell_1\br u_2\over\to
\ell_2\br u_1\over\to X$, where $u_1$ and $u_2$ are compact operators. Let
now $Z$ be any Banach 
space. By a result of S.Bellenot [1] (J.S.Morell and J.R.Retherford [12])
there exists a quotient
space $Q$ of $Z^*$ and compact operators $v_1:Q\to X$ and $v_2:\ell_2\to Q$
such that $u_1=v_1
v_2$. Invoke the compact lifting property of $\ell_1$ to finish the proof. 
\ms
The situation resembles the one encountered in Proposition 2. Question: Can
one get a 
factorization through $Z$ rather than through $Z^*$?

\vglue 10mm\n
\centerline {\SUBT REFERENCES}
\bs\sm
{\aop
\item{[1]~} S.Bellenot:\quad The Schwartz-Hilbert variety.\par 
Mich.Math.J. {\bf 22} (1975) 373~-377.
\sm
\item{[2]~} F.F.Bonsall and J.Duncan:\quad Complete Normed Algebras.\par 
Springer-Verlag 1973.
\sm
\item{[3]~} C.Cho and W.B.Johnson:\quad A characterization of subspaces $X$
of $\ell^p$ for which $K(X)$ is an M-ideal in $L(X)$. Proc.AMS {\bf 93}
(1985) 466-470.
\sm
\item{[4]~} P.J.Cohen:\quad Factorization in group algebras.\par 
Duke Math.J. {\bf 26} (1959) 199-205.
\sm
\item{[5]~} H.G.Dales and H.Jarchow:\quad Continuity of homomorphisms and
derivations from algebras of approximable and nuclear operators.\par
Math.Proc.Cambr.Phil.Soc., to appear. 
\sm
\item{[6]~} P.G.Dixon:\quad Left approximate identities in algebras of
compact operators
on Banach spaces. Proc.Roy.Soc.Edinburgh {\bf 104A} (1989)169~-175.
\sm
\item{[7]~} N.Gr\oe nb\ae k and G.A.Willis:\quad Approximate identities in
the Banach algebras
of compact operators. Canad.Math.Bull.
\sm
\item{[8]~} A.~Grothendieck: Produits tensoriels topologiques at espaces
nucl\'eaires.\par
Mem.Amer.Math.Soc. {\bf 16}, 1955.
\sm
\item{[9]~} K.John:\quad On the compact non-nuclear problem.\par 
Math.Ann. {\bf 287} (1990) 509-514.
\sm
\item{[10]~} W.B.Johnson:\quad A complementary universal conjugate Banach
space and its relation 
to the approximation problem.\par 
Israel J.Math. {\bf 13} (1972) 301-310. 
\sm
\item{[11]~} J.Lindenstrauss and L.Tzafriri; Classical Banach Spaces I,
Sequence Spaces.\par
Springer-Verlag, Erg.Math.Grenzg. {\bf 92}, 1977. 
\sm
\item{[12]~} J.S.Morell and J.R.Retherford:\quad p-trivial Banach spaces.\par
Studia Math. {\bf 43} (1972) 1-25.
\sm
\item{[13]~} A.Pietsch:\quad Operator Ideals.\par 
VEB Deutscher Verlag der Wissenschaften 1978 / North-Holland 1980. 
\sm
\item{[14]~} G.Pisier:\quad Un th\'eor\`eme sur les op\'erateurs entre
espaces de Banach qui se
factorisent par un espace de Hilbert.\par 
Ann. \'Ecole Norm.Sup. {\bf 13} (1980) 23-43.
\sm
\item{[15]~} G.Pisier:\quad Counterexample to a conjecture of Grothendieck.\par
Acta Math. {\bf 151} (1983) 180-208.
\sm
\item{[16]~} C.Samuel:\quad Bounded approximate identities in the algebra
of compact operators 
on a Banach space.\par Proc.AMS {\bf 117} (1993) 1093-1096. 
\sm
\item{[17]~} Y.V.Selivanov:\quad Homological characterizations of the
approximation property 
for Banach spaces.\par 
Glasgow Math.J. {\bf 34} (1992) 229-239. 
\sm
\item{[18]~} A.Szankowski:\quad Subspaces without the approximation
property.\par
Israel J.Math. {\bf 30} (1978) 123-129.
\sm
\item{[19]~} G.A.Willis:\quad The compact approximation property does not
imply the approximation
property.\par 
Studia Math. {\bf 103} (1992) 99-108.

\par}

\vskip 10mm\n
\halign to\hsize{#\hfil&#\hskip 40mm&#\hfil\cr P.G.Casazza&&H.Jarchow\cr
Department of Mathematics&&Mathematisches Institut\cr 
University of Missouri-Columbia&&Universit\"at Z\"urich\cr 
Columbia, Mo 65211&&CH 8057 Z\"urich\cr
USA&&Switzerland\cr
e-mail:&&e-mail:\cr
Pete@casazza.cs.missouri.edu&&Jarchow@math.unizh.ch\cr}

\bye